\documentclass[review]{amsart}
\usepackage{graphicx}
\newtheorem{thm}{Theorem}[section]
\newtheorem{cor}[thm]{Corollary}
\newtheorem{lem}[thm]{Lemma}
\newtheorem{prop}[thm]{Proposition}
\theoremstyle{definition}

\theoremstyle{remark}

\numberwithin{equation}{section}
\newcommand{\Real}{\mathbb R}

\def\R{{\bf R}}

\def\N{{\bf N}}

\def\d{\displaystyle}

\def\p{\partial}

\begin{document}

\title[Supercritical Dissipative Wave Equations]
{Global Regularity for Supercritical Nonlinear Dissipative Wave Equations in 3D}
\author{Kyouhei Wakasa and Borislav Yordanov}

\address{\begin{tabular}{ccc} Muroran Institute of Technology& $\quad$ $\quad$
                                                           & Hokkaido University, Sapporo    \\
                                                           &  &  and                            \\
                                                           &  & Institute of Mathematics, Sofia \\
                                                           &  &
\end{tabular}}

\email{wakasa@mmm.muroran-it.ac.jp\qquad \quad byordanov@oia.hokudai.ac.jp}

%\thanks{}%
%\subjclass{}%
%\keywords{}%

%\date{}%
%\dedicatory{}%
%\commby{*}%

\date{\today}
\subjclass{Primary 35L70; Secondary 35A05} \keywords{wave
equation, nonlinear damping, supercritical, regularity}

% ----------------------------------------------------------------
\begin{abstract}

The nonlinear wave equation $u_{tt}-\Delta u +|u_t|^{p-1}u_t=0$ is
shown to be globally well-posed in the Sobolev spaces of radially symmetric functions $H^k_{\rm
rad}(\Real^3)\times H^{k-1}_{\rm rad}(\Real^3)$ for all $p\geq 3$ and $k\geq 3$.
Moreover, global $C^\infty $ solutions are obtained when
the initial data are $C_0^\infty$ and exponent $p$ is an odd integer.

The radial symmetry allows a reduction to the one-dimensional case where an important observation of A. Haraux~\cite{Ha}
can be applied, i.e., dissipative nonlinear wave equations contract initial data in $W^{k,q}(\Real)\times W^{k-1,q}(\Real)$ for
all $k\in[1,2]$ and $q\in [1,\infty]$.

\end{abstract}
\maketitle

% ----------------------------------------------------------------
\section{Introduction}

Dissipative nonlinear wave equations are well-posed in $H^k(\Real^n)\times H^{k-1}(\Real^n)$
for all $k\in [1,2]$ and $n\geq 1$ under the monotonicity condition of Lions and Strauss~\cite{LiSt}.
This global result makes no exception for nonlinear dissipations of supercritical power, as determined by invariant
scaling, and gets around the stringent conditions for well-posedness of general nonlinear wave equations; see
Ponce and Sideris~\cite{PS}, Lindblad~\cite{Lind}, Wang and Fang~\cite{WaFa} and the references therein.

The monotonicity method has been less effective in studying higher regularity.
It is still an open question whether
supercritical problems are globally well-posed in Sobolev spaces with index $k>2$.
The purpose of this paper is to give an affirmative answer when $n=3$ and initial data have radial symmetry.
To state the result, let $\Box u =u_{tt}-\Delta u$ be the d'Alembertian in $\Real^{3+1}$ and $Du=(\nabla u,\partial_t)$
be the space-time gradient of $u$.  Standard notations are also $\|\cdot\|_q$, for the norm in $L^q(\Real^3)$ with $q\in[1,\infty]$,
and $D^\alpha$, for the partial derivative of integer order $\alpha=(\alpha_1,\alpha_2,\alpha_3,\alpha_4)$.

We consider the dissipative nonlinear wave equation
\begin{eqnarray}
\label{main}
 \Box u+|u_t|^{p-1}u_t=0,& & x\in \Real^3,\quad t>0,
 \end{eqnarray}
with the Cauchy data
\begin{eqnarray}
\label{data}
 u|_{t=0}=u_0,\qquad u_t|_{t=0}=u_1,& & x\in \Real^3.
 \end{eqnarray}
Our main assumptions are $p\geq 3$ and $(u_0,u_1)\in H^{k}_{\rm rad}(\Real^3)\times H^{k-1}_{\rm rad}(\Real^3)$, where
the radially symmetric Sobolev spaces are defined as
\begin{eqnarray*}
H^k_{\rm rad}(\Real^3)=\{u \in H^k(\Real^3): (x_i\partial_{x_j}-x_j\partial_{x_i})u(x)=0\ \hbox{ for } \ 1\leq i<j\leq 3\}.
\end{eqnarray*}
Clearly, such $u(x)$ depend only on $|x|=(x_1^2+x_2^2+x_3^2)^{1/2}$. The radially symmetric spaces are invariant
under the evolution determined by (\ref{main}), (\ref{data}); see \cite{St}, \cite{ShSt2}.

\begin{thm}
 \label{th1}
Assume that $p\geq 3$ and $(u_0,u_1)\in H^3_{\rm
rad}(\Real^3)\times H^2_{\rm rad}(\Real^3).$ Then problem
(\ref{main}), (\ref{data}) admits a unique global solution $u,$ such
that
\[
D^\alpha u \in C([0,\infty), H^{3-|\alpha|}_{\rm rad}(\Real^3)),\qquad
 |\alpha|\leq 3.
\]
\end{thm}

\begin{cor}
 \label{cor011}
The global solution of problem (\ref{main}), (\ref{data}), given by Theorem~\ref{th1}, satisfies
the following uniform estimates: for $t\geq 0$,
\begin{eqnarray*}
 \sum_{1\leq |\alpha|\leq 2}\int_0^t \|D^\alpha u(s)\|^2_{\infty}\: ds
  & \leq & C_1(u_0,u_1), \\
 \sum_{|\alpha|\leq 2}\|DD^\alpha u(t)\|_2 & \leq  & C_2(u_0, u_1),
\end{eqnarray*}
where $C_j(u_0,u_1),$ $j=1, 2,$ are finite whenever $\|u_0\|_{H^{3}}+\|u_1\|_{H^{2}}$ is finite.
\end{cor}

It is easy to derive the propagation of higher regularity from Theorem~\ref{th1} and Sobolev embedding inequalities,
if the nonlinearity allows further differentiation.

\begin{thm}
\label{th01}
Let $p$ be an odd integer and $(u_0,u_1)\in H^k_{\rm
rad}(\Real^3)\times H^{k-1}_{\rm rad}(\Real^3)$ with an integer $k\geq 4$.
Problem (\ref{main}), (\ref{data}) admits a unique global solution $u,$ such
that
\[
D^\alpha u \in C([0,\infty), H^{k-|\alpha|}_{\rm rad}(\Real^3)),\quad
|\alpha|\leq k.
\]
In addition, the following estimates hold uniformly in $t\geq 0$:
\begin{eqnarray*}
 \sum_{1\leq |\alpha|\leq k-1}\int_0^t \|D^\alpha u(s)\|^2_{\infty}\: ds
  & \leq & C_1^{(k)}(u_0,u_1), \\
 \sum_{|\alpha|\leq k-1}\|DD^\alpha u(t)\|_2 & \leq  & C_2^{(k)}(u_0, u_1),
\end{eqnarray*}
where $C_j^{(k)}(u_0,u_1),$ $j=1, 2,$ are finite whenever $\|u_0\|_{H^{k}}+\|u_1\|_{H^{k-1}}$ is finite.

If $(u_0,u_1)\in C^\infty_{\rm rad}(\Real^3)\times C^\infty_{\rm
rad}(\Real^3)$ and $(u_0,u_1)$ are compactly supported, then problem
(\ref{main}), (\ref{data}) admits a unique global solution
\[
u \in C^\infty([0,\infty), C^\infty_{\rm rad}(\Real^3)),
\]
such that $u(\cdot,t)$ is also compactly support for all $t\geq 0$.
\end{thm}

The above results rely heavily on the dissipative nonlinearity and radial symmetry.
In less favorable circumstances, global regularity is well beyond the reach of existing methods which
require the nonlinearity be slightly weaker than the linear differential operator.
Nonlinear conservative wave equations and Schrodinger equations, for example, are well-understood only in
subcritical and critical cases; see Struwe~\cite{StM}, Grillakis~\cite{Gr1}, Shatah and Struwe~\cite{ShSt1},
Nakanishi~\cite{Na} for the former and Bourgain~\cite{Bo}, Grillakis~\cite{Gr2}, Colliander, Keel, Staffilani, Takaoka
and Tao~\cite{CoKeStTaTa} for the latter.
The supercritical nonlinearities allow mostly conditional results about the equivalence of certain norms and nature of blow up, such as Kenig and Merle~\cite{KM} and Killip and Visan \cite{KV}.
Exceptions are the well-posedness for log-supercritical nonlinearities and radial data by Tao~\cite{Ta3}
and the well-posdeness and scattering for log-log-supercritical nonlinearities by Roy \cite{Roy2}.

Dissipative and conservative nonlinear wave equations behave very differently
even in Sobolev spaces with $k\leq 2$ derivatives. This is evident from the invariant scaling of (\ref{main}):
$u(x,t)\mapsto u_\lambda(x,t)=\lambda^{(2-p)/(p-1)}u(\lambda x,\lambda t)$, with $\lambda>0$. Then
\begin{equation}
\label{scaling}
\|u_\lambda(t)\|_{\dot{H}^k}=\lambda^{(2-p)/(p-1)+k-3/2}\|u(\lambda t)\|_{\dot{H}^k},
\end{equation}
so the critical space for well-posedness is $\dot{H}^{k_c}(\Real^3)\times \dot{H}^{k_c-1}(\Real^3)$ with
\[
k_c=3/2+(p-2)/(p-1).
\]
However, the result of \cite{LiSt} shows that $k\in [1,2]$ is sufficient for any $p>1$.
Here the monotone dissipation plays a decisive role, since (\ref{main}) remains dissipative after applying $D^\alpha$ with $|\alpha|=1$:
\begin{equation}
\label{diff.main}
\Box D^\alpha u+p|u_t|^{p-1}D^\alpha u_t=0.
\end{equation}
We can use $u_t$ and $D^\alpha u_t$ as multipliers for (\ref{main}) and (\ref{diff.main}), respectively, to obtain
\begin{eqnarray*}
  \|Du(t)\|^2_{2} + 2\int_0^t \|u_s(s)\|^{p+1}_{p+1}\: ds & = & \|Du(0)\|^2_{2},\\
  \|DD^\alpha u(t)\|^2_{2} + 2p\int_0^t \||u_s(s)|^{(p-1)/2}D^\alpha u_s(s)\|^2_{2}\: ds & = & \|DD^\alpha u(0)\|^2_{2}.
\end{eqnarray*}
As derivatives of order $k\leq 2$ turn out to be {\em a priori}
bounded, the corresponding well-posedness results readily follow from the monotonicity method
of \cite{LiSt}.

Two differentiations of (\ref{main}) yield an equation that is no longer dissipative.
It becomes a nontrivial task to derive uniform estimates for derivatives of order three and higher.
Some information is provided by the invariant scaling, which predicts the non-concentration
of second-order norms if $k=2$ and $(2-p)/(p-1)+1/2>0$ in (\ref{scaling}). Thus, the range of subcritical exponents is $p<3$.
For the complementary range $p\geq 3$ and $k>2$, the global well-posedness of (\ref{main})
in $H^k(\Real^3)\times H^{k-1}(\Real^3)$ is a difficult problem.
Until recently, the proof has been known only in the critical case $p=3$ with radially symmetric data~\cite{TUY}.

This paper shows that no critical exponent exists for the regularity of dissipative nonlinear wave equations with radial symmetry.
The development of singularities is prevented by the monotonicity of seminorms involving
second-order derivatives. We actually obtain, after setting $r=|x|$ and
$X_{\pm}(r,t)=(\partial_t \pm \partial_r)(r\p_tu)$, that
\[
\max_{\pm}\sup_{r>0}|X_{\pm}(r,t)| \leq  \max_{\pm}\sup_{r>0}|X_{\pm}(r,0)|, \quad t\geq 0.
\]
Such decreasing quantities exist only for linear and dissipative nonlinear wave equations
in $\Real^3$ with radial data. The proof is based on differentiating (\ref{main}) in $t$ and rewriting the equation for $\partial_tu$ as
\[
(\partial_t \mp \partial_r)X_{\pm}+\frac{p}{2}|\partial_t u|^{p-1}(X_++X_-)=0.
\]
A. Haraux \cite{Ha} has applied the same idea to derive $W^{2,\infty}(\Real)\times W^{1,\infty}(\Real)$ estimates in the one-dimensional case.
In higher dimensions with radial symmetry, similar $W^{1,\infty}(\Real)\times L^{\infty}(\Real)$ estimates away from $r=0$ are used by
Joly, Metivier and Rauch~\cite{JoMeRa}, Liang~\cite{Li} and Carles and Rauch~\cite{CaRa}.
Since the one-dimensional reduction does not work for general data, the global well-posedness in $H^{k}(\Real^3)\times H^{k-1}(\Real^3)$, with $k>2,$
remains completely open.

The rest of this paper is organized as follows. Section~2 contains
several basic facts and estimates for the wave equation with radially symmetric data in
$\Real^3$. The regularity problem in $H^3_{\rm rad}(\Real^3)\times H^{2}_{\rm rad}(\Real^3)$
is studied in Section~3.
In the final Section~4, we establish Theorem~1.3 about well-posedness in $H^k_{\rm
rad}(\Real^3)\times H^{k-1}_{\rm rad}(\Real^3)$ with $k\geq 4$ and
$C^\infty_{\rm rad}(\Real^3)\times C^\infty_{\rm rad}(\Real^3)$ with
compact support.

\section{Basic Estimates}

First of all, we state the energy estimates and
Strichartz estimates for the radial wave equation in $\Real^3\times \Real$.

\begin{lem}
\label{lem02} Let $u$ be a solution of the Cauchy problem
 in $\Real^3\times \Real$
\[
\Box u=F,\qquad u|_{t=0}=u_0,\qquad u_t|_{t=0}=u_1.
\]

(a) For any source and initial data, $u$ satisfies the energy
estimate
\begin{eqnarray*}
\|Du(t)\|_{2} \leq C(\|\nabla u_0\|_{2}+\|u_1\|_{2})+C\int_0^t
\|F(s)\|_{2}\: ds
\end{eqnarray*}
with an absolute constant $C$ for all $t\geq 0.$

(b) For radial source and initial data, $u$ is also a radial
function which satisfies
\begin{eqnarray*}
\left(\int_0^t \|u(s)\|_{\infty}^2 \: ds \right)^{1/2} \leq
C(\|\nabla u_0\|_{2}+\|u_1\|_{2})+C\int_0^t \|F(s)\|_2\: ds
\end{eqnarray*}
for all $t\geq 0.$
\end{lem}

Part (a) can be found in Strauss~\cite{St}, H\"{o}rmander~\cite{Ho}, and Shatah and Struwe~\cite{ShSt2}.
Estimate (b) is the so-called ``radial Strichartz estimate" in $3D$.
Klainerman and Machedon~\cite{KlMa} have found the
homogeneous version of (b) which implies the non- homogeneous estimate stated here.

The following is a collection of useful facts concerning local solvability and other
properties of problem (\ref{main}), (\ref{data}).
These results can be found in \cite{St}, \cite{Ho} and \cite{ShSt2}.

\begin{lem}
\label{lem00} Let $k\geq 3$ and
 $(u_0,u_1)\in H^k(\Real^3)\times H^{k-1}(\Real^3).$

(a) There exists $T>0,$ such that problem (\ref{main}), (\ref{data})
has a unique solution $u$ satisfying
\[
D^\alpha u \in C([0,T], H^{k-|\alpha|}(\Real^3)),\qquad |\alpha|\leq
k.
\]
Moreover, we have
\[
\sup_{t\in [0,T]}\|D^\alpha u(t)\|_2\leq C_k,
\]
where $T$ and $C_k$ can be chosen to depend continuously on
$\|u_0\|_{H^k}+\|u_1\|_{H^{k-1}}.$

(b) The continuation principle holds: if $T_*=T_*(u_0,u_1)$ is the
supremum of all numbers $T$ for which (a) holds, then either
$T_*=\infty$ or
\[
\sup_{t\in [0,T_*)}\|D^\alpha u(t)\|_2=\infty
\]
for some $\alpha$ with $|\alpha|\leq k.$

(c) If the data $(u_0,u_1)$ are radially symmetric, the solution
$(u,u_t)$ is also radially symmetric.

\end{lem}

Next, we state two preliminary estimates for problem (\ref{main}),
(\ref{data}). These results, called the energy dissipation laws, are already discussed in the
introduction.
\begin{lem}
\label{lem01} Assume that
 $(u_0,u_1)\in H^3(\Real^3)\times H^{2}(\Real^3)$ and $|\alpha|=1.$
Let $u$ be the solution of problem (\ref{main}), (\ref{data})
for $t\in [0,T]$, given by Lemma~\ref{lem00}. Then
\begin{eqnarray*}
\frac{1}{2}\|Du(t)\|_2^2+\int_0^t \|u_s(s)\|^{p+1}_{p+1}\: ds & = & \frac{1}{2}\|Du(0)\|_2^2,\\
\frac{1}{2}\|DD^\alpha u(t)\|_2^2+p\int_0^t \||u_s(s)|^{(p-1)/2}D^\alpha
u_s(s)\|^2_{2}\: ds & = & \frac{1}{2}\|DD^\alpha u(0)\|_2^2,
\end{eqnarray*}
for all $t\in[0,T].$ Thus, the following norms of $u$ are uniformly bounded:
\begin{eqnarray*}
& & \|Du(t)\|_2\leq \|Du(0)\|_2,\quad \ \|DD^\alpha u(t)\|_2\leq
\|DD^\alpha u(0)\|_2,\quad |\alpha|\leq 1,\\
& & \int_0^T(\|u_s(s)\|^{p+1}_{p+1}+\||u_s(s)|^{(p-1)/2}D^\alpha u_s(s)\|^{2}_{2})\: ds\leq \|Du(0)\|_2^2+\|DD^\alpha u(0)\|_2^2.
\end{eqnarray*}
\end{lem}

{\bf Proof.}
Multiplying equation (\ref{main}) by $u_t$, we get
\begin{eqnarray*}
%\label{divE}
0=(\Box u+|u_t|^{p-1}u_t)u_t = \left(\frac{|Du|^2}{2}\right)_t
-\hbox{div}(u_t\nabla u) +|u_t|^{p+1}.
\end{eqnarray*}
The first-order energy identity follows from the integration on $\Real^3$ and divergence theorem
if $u(x,t)$ has compact support with respect to $x.$
More generally, we can
approximate $(u_0(x),u_1(x))$ with compactly supported $C^\infty$
functions and use the finite propagation speed to show that the
boundary integral of $\hbox{div}(u_t\nabla u)$ is zero. Property (a)
in Lemma~\ref{lem00} implies that the approximations will converge
to the actual solution.

Similarly, we can differentiate equation (\ref{main}) and multiply
with $D^\alpha u_t(x,t)$ to show the second-order energy identity. Let us recall
that we work with solutions whose third-order derivatives belong to
$L^2(\Real^3).$
 \qed

Finally, we state the strong version of Strauss inequality. The original version
can be found as Radial Lemma 1 of Strauss \cite{St1}.

\begin{lem}
\label{lem07}
(Strong version of Strauss Inequality \cite{St1}) Let $U\in H^{1}_{\rm rad}(\Real^3)$. There exists a constant $C>0$,
such that for every $R>0$
\[
R |U(R)|\le C\|U\|_{H^1(|x|>R)}\rightarrow 0 \ \hbox{as}\ R\rightarrow\infty.
\]
\end{lem}

\section{Global existence of radial solutions in $H^{3}\times H^2$}

The following lemma is essential for the proof of Theorem~\ref{th1}.
This approach to $L^\infty$ estimates of second-order derivatives is borrowed from Haraux~\cite{Ha}.

\begin{lem}
\label{lem:Haraux} Let $(u_0,u_1)\in H^3_{\rm rad}(\Real^3)\times H^{2}_{\rm rad}(\Real^3)$ and
let $u$ be the local solution of problem (\ref{main}), (\ref{data})
constructed in  Lemma \ref{lem00}. There exists a positive constant
$C=C(\|u_0\|_{H^3},\|u_1\|_{H^2})$,
such that
\begin{eqnarray}
& & \sup_{r\in[0,\infty)}r|\p_t^2u(r,t)|\le C,\label{H_Est1}\\
& & \sup_{r\in[0,\infty)} |\p_t u(r,t)|\le C,\quad
\sup_{r\in[0,\infty)}r|\p_r\p_{t}u(r,t)|\le C\label{H_Est2},
\end{eqnarray}
for all $t\in[0,T_*)$.
\end{lem}
{\bf Proof.} Since the solution of (\ref{main}) is radially symmetric, the equation becomes
\begin{eqnarray}
\label{main:rad}
(\p^2_{t}-\partial^2_{r})(ru)+r|\p_tu|^{p-1}\p_tu=0,
\end{eqnarray}
where $r=|x|$. Let us define
\[
X_{\pm}(r,t)=(\partial_t \pm \partial_r)(r\p_tu).
\]
We differentiate (\ref{main:rad}) in $t$ to obtain an equation for $\partial_t u$.
Then we multiply through by $X_{\pm}^{2k-1}$ with $k\in\N$ to obtain
\[
\frac{(\partial_t-\partial_r)}{2k}\{X_+(r,t)\}^{2k}
+pr|\partial_tu|^{p-1}\partial^2_t u\{X_{+}(r,t)\}^{2k-1}=0
\]
and
\[
\frac{(\partial_t+\partial_r)}{2k}\{X_-(r,t)\}^{2k}
+pr|\partial_tu|^{p-1}\partial^2_tu\{X_{-}(r,t)\}^{2k-1}=0.
\]
Adding the above two equations and using $X_+(r,t)+X_-(r,t)=2r\partial^2_tu$, we have
\[
\frac{(\partial_t-\partial_r)}{2k}X_+^{2k}+\frac{(\partial_t+\partial_r)}{2k}X_-^{2k}
+p|\partial_tu|^{p-1}\frac{X_++X_-}{2}(X_+^{2k-1}+X_-^{2k-1})=0.
\]
Here we remark that the following inequality holds for $a,b\in \R$ and $k \in \N$:
\[
(a+b)(a^{2k-1}+b^{2k-1})\ge0.
\]
From this inequality with $a=X_+$ and $b=X_-$, we derive
\begin{eqnarray*}
\partial_t(X_+^{2k}+X_-^{2k})+\partial_r(X_-^{2k}-X_+^{2k})\le 0.
\end{eqnarray*}
It is easy to see that integration on $[0,R]\times[0,t]$ implies
\begin{eqnarray*}
\int_{0}^{R}\{X_+^{2k}(r,t)+X_-^{2k}(r,t)\}dr & \le &\int_{0}^{R}\{X_+^{2k}(r,0)+X_-^{2k}(r,0)\}dr\\
& & +\int_{0}^{t}\{X_+^{2k}(R,s)-X_-^{2k}(R,s)\}ds,
\end{eqnarray*}
since $X_-^{2k}(0,t)-X_+^{2k}(0,t)=0$. We also notice that $|X_{\pm}(R,t)|\rightarrow 0$ as
$R\rightarrow \infty$, which is a consequence of Lemma \ref{lem07}. Thus, we get
\[
\|X_{+}(\cdot,t)\|_{L^{2k}(\Real_+)}^{2k}+\|X_{-}(\cdot,t)\|_{L^{2k}(\Real_+)}^{2k}
\le\|X_{+}(\cdot,0)\|_{L^{2k}(\Real_+)}^{2k}+\|X_{-}(\cdot,0)\|_{L^{2k}(\Real_+)}^{2k}.
\]
The two terms on the right hand side satisfy
\[
\|X_{\pm}(\cdot,0)\|_{L^{2k}(\Real_+)}^{2k}\leq \|X_{\pm}(\cdot,0)\|_{L^{\infty}(\Real_+)}^{2k-2}\|X_{\pm}(\cdot,0)\|_{L^2(\Real_+)}^2,
\]
where $\|X_{\pm}(\cdot,0)\|_{L^2(\Real_+)}<\infty$ by Lemma~\ref{lem07}
and $(u_0,u_1)\in H^3_{\rm rad}(\Real^3)\times H^2_{\rm rad}(\Real^3)$.
This observation yields a more convenient estimate:
\begin{eqnarray*}
\|X_{+}(\cdot,t)\|_{L^{2k}(\Real_+)}^{2k}+\|X_{-}(\cdot,t)\|_{L^{2k}(\Real_+)}^{2k}&\le& \|X_+(\cdot,0)\|_{L^{\infty}(\Real_+)}^{2k-2}\|X_+(\cdot,0)\|_{L^2(\Real_+)}^2\\
&&+\|X_-(\cdot,0)\|_{L^{\infty}(\Real_+)}^{2k-2}\|X_-(\cdot,0)\|_{L^2(\Real_+)}^2.
\end{eqnarray*}
Letting $k\rightarrow\infty$, we have that
\begin{eqnarray}
\label{boundedness}
& & \max\{\|X_+(\cdot,t)\|_{L^{\infty}(\Real_+)},\|X_-(\cdot,t)\|_{L^{\infty}(\Real_+)}\}\nonumber\\
& & \le \max\{\|X_+(\cdot,0)\|_{L^{\infty}(\Real_+)},\|X_-(\cdot,0)\|_{L^{\infty}(\Real_+)}\}.
\end{eqnarray}

The remaining part of the proof is standard.
Since $X_{+}(r,t)+X_{-}(r,t)=2r\p^2_{t}u$, claim (\ref{H_Est1})
follows from (\ref{boundedness}).
We can write $2\p_r(r\p_{t}u)=X_{+}(r,t)-X_{-}(r,t)$ and
\[
2r|\p_tu(r,t)|\le \int_{0}^{r}|X_{+}(\rho,t)-X_{-}(\rho,t)|d\rho\le Cr.
\]
Hence $\|\p_tu(\cdot,t)\|_{L^{\infty}(\Real_+)}\leq C/2<\infty$,
which is the first claim in (\ref{H_Est2}).
We finally observe that $2r\p_r\p_t u(r,t)=X_{+}(r,t)-X_{-}(r,t)-2\p_tu(r,t)$, so the second claim in
(\ref{H_Est2}) follows from the first one and (\ref{boundedness}). The proof is complete.
\qed

\vskip 0.5truecm

{\bf Proof of Theorem~\ref{th1}.} It is sufficient to show that $\|u(t)\|_{H^3}+\|u_t(t)\|_{H^2}$
does not blow up as $t\rightarrow T_*$, where $u$ is the local solution of (\ref{main}) given by
Lemma~\ref{lem00}. The first and second order norms of $u$ are {\it a priori} bounded from Lemma~\ref{lem01},
so the global existence of $u$ is guaranteed by the next result about third order norms.

\begin{prop}
\label{pr05} Assume that
 $(u_0,u_1)\in H^3_{\rm rad}(\Real^3)\times H^{2}_{\rm rad}(\Real^3)$ and
let $u$ be the local solution of problem (\ref{main}), (\ref{data})
constructed in  Lemma~\ref{lem00}. There exists a positive constant $C=C(u_0,u_1,T_*)$
such that
\begin{eqnarray*}
\sum_{|\alpha|=3}\|D^{\alpha}u(t)\|_2  \le C(u_0,u_1,T_*)
\end{eqnarray*}
for all $t\in [0,T_*).$
\end{prop}

 {\bf Proof.} Differentiating (\ref{main}) twice, we find that
 $D^\alpha u$ is a weak solution of
\[
 \Box D^\alpha u+p|u_t|^{p-1}D^\alpha u_t+p(p-1)|u_t|^{p-3}u_tD^{\alpha_1}u_tD^{\alpha_2}u_t=0,
\]
where $\alpha=\alpha_1+\alpha_2$ with $|\alpha_1|=1$ and $|\alpha_2|=1$.
Thus, $D^\alpha u$ satisfies the energy
estimate in Lemma~\ref{lem02} (a):
\begin{eqnarray*}
 \|DD^\alpha u(t)\|_2 & \leq & C\|DD^\alpha u(0)\|_2\\
                      &      & + C\int_0^t \|u_s^{p-1}(s)D^\alpha u_s(s)\|_2 \: ds\\
                      &      &+C\int_0^t \|u_s^{p-2}(s) D^{\alpha_1}u_s(s)D^{\alpha_2}u_s(s)\|_2\: ds
\end{eqnarray*}
for $t\in [0,T_*)$. We notice that
\[
\|u_s^{p-1}(s)D^\alpha u_s(s)\|_2\le C \|D^{\alpha}u_s(s)\|_2,
\]
by (\ref{H_Est1}) in Lemma~\ref{lem:Haraux}, and
\[
\|u_s^{p-2}(s) D^{\alpha_1}u_s(s)D^{\alpha_2}u_s(s)\|_2\le C \left\|\frac{D^{\alpha_1}u_s(s)}{|\cdot|}\right\|_2,
\]
by (\ref{H_Est1}) and (\ref{H_Est2}) in Lemma~\ref{lem:Haraux}. Thus,
\begin{eqnarray*}
 \|DD^\alpha u(t)\|_2 & \leq & C\|DD^\alpha u(0)\|_2\\
                      &      & + C\int_0^t \|D^\alpha u_s(s)\|_2 \: ds\\
                      &      &+C\int_0^t \left\|\frac{D^{\alpha_1}u_s(s)}{|\cdot|}\right\|_2 \: ds.
\end{eqnarray*}
The third term of the right hand can be estimated by Hardy's inequality:
\begin{eqnarray*}
 \|DD^\alpha u(t)\|_2 & \leq & C\|DD^\alpha u(0)\|_2\\
                      &      & + C\int_0^t \|D^\alpha u_s(s)\|_2 \: ds\\
                      &      &+C\int_0^t\|DD^{\alpha_1}u_s(s)\|_2 \: ds.
\end{eqnarray*}
We add these estimates for all $|\alpha|=2$ to get
\[
\sum_{|\alpha|=3}\|D^{\alpha}u(t)\|_2 \:  \leq C\sum_{|\alpha|=3}\|D^{\alpha}u(0)\|
+C\int_{0}^{t}\sum_{|\alpha|=3}\|D^\alpha u(s)\|_2ds.
\]
Making use of Gronwall's inequality, we finally have
\[
\sum_{|\alpha|=3}\|D^{\alpha}u(t)\|_2  \le C\sum_{|\alpha|=3}\|D^{\alpha}u(0)\|
e^{Ct}.
\]
This completes the proof of global existence.
\qed

{\bf Proof of Corollary~\ref{cor011}.} We can now verify the uniform estimates of $L^2$ norms and square integrability of
$L^{\infty}$ norms on $[0,\infty)$.
Notice that inequality (b) in Lemma~\ref{lem02} holds on any interval $[t_0,t]$ for all $|\alpha|=1$.
Hence we get
\begin{eqnarray*}
  \sum_{|\alpha|=1}\left(\int_{t_0}^t \|D^\alpha u(s)\|_\infty^2\: ds\right)^{1/2}
   & \leq & C\sum_{|\alpha|=1}\|DD^\alpha u(t_0)\|_2\\
   &      &  +C\sum_{|\alpha|=1}\int_{t_0}^t \|u_s^{p-1}(s)D^\alpha u_s(s)\|_2\: ds.
\end{eqnarray*}
It is convenient to abbreviate the left hand side as
\[
N_1(t)=\sum_{|\alpha|=1}\left(\int_{t_0}^t \|D^\alpha u(s)\|_\infty^2\: ds\right)^{1/2},\qquad t\geq t_0,
\]
and estimate the integrand on the right hand aside as
$$\|u_s^{p-1}(s)D^\alpha u_s(s)\|_2\le C\|D^{\alpha}u(s)\|_{\infty}\|u_s^{(p-1)/2}(s)D^\alpha u_s(s)\|_2,$$
from (\ref{H_Est2}). Applying the Cauchy inequality on $[t_0,t]$ to the initial estimate,
\begin{eqnarray*}
N_1(t) & \leq & C\sum_{|\alpha|=1}\|DD^\alpha u(t_0)\|_2\\
   &      &  +CN_1(t)\left(\int_{t_0}^t \|u_s^{(p-1)/2}(s)Du_s(s)\|_2^2\: ds\right)^{1/2}.
\end{eqnarray*}
The above integral converges on $[0,\infty)$ by Lemma~\ref{lem01}. We can find $t_0,$ such that
\begin{equation}
\label{conv}
 C\left(\int_{t_0}^t \|u_s^{(p-1)/2}(s)Du_s(s)\|_2^2\:
ds\right)^{1/2}\leq \frac{1}{2},\qquad t\geq t_0.
\end{equation}
Thus, $\d N_1(t)\leq 2C\sum_{|\alpha|= 1}\|DD^\alpha u(t_0)\|_2$ for
all $t\geq t_0.$ If $t<t_0$, a similar estimate for $N_1(t)$ follows from
Theorem~\ref{th1}. Hence $\|D^\alpha u(\cdot)\|_\infty\in L^2([0,\infty))$ for $|\alpha|=1$.

Next, we show that the second order norms in Corollary~\ref{cor011} are also uniformly bounded.
Set
\[
N_2(t)=\sum_{|\alpha|= 2}\left(\int_{t_0}^t \|D^\alpha
u(s)\|_\infty^2\: ds\right)^{1/2},\qquad t\geq t_0,
\]
and apply estimate (b) in Lemma \ref{lem02} to obtain
\begin{eqnarray*}
 \left(\int_{t_0}^t \|D^\alpha u(s)\|_{\infty}^2 \: ds
\right)^{1/2} & \leq & C\|DD^\alpha u(t_0)\|_{2}+C\int_{t_0}^t\|u_s^{p-1}(s)D^\alpha u_s(s)\|_2\: ds\\
 & & +C\int_{t_0}^t\|u_s^{p-2}(s)(Du_s(s))^2\|_2\: ds,
\end{eqnarray*}
where $|\alpha|=2.$ Noticing that (\ref{H_Est2}) gives
\begin{eqnarray}
\|u_s^{p-1}(s)D^\alpha u_s(s)\|_2&\le &C\|u_s(s)\|_{\infty}^2\|D^{\alpha}u_s(s)\|_{2},\label{uni_est1}\\
\|u_s^{p-2}(s)(Du_s(s))^2\|_2 &\le &C\|D^{\alpha}u(s)\|_{\infty}\|u_s^{(p-1)/2}(s)Du_s(s)\|_2\label{uni_est2}
\end{eqnarray}
and applying the Cauchy inequality on $[t_0,t]$, we get
\begin{eqnarray*}
N_2(t) & \leq & C\sum_{|\alpha|=2}\|DD^\alpha u(t_0)\|_2\\
       &      & +C\int_{t_0}^t\|u_s(s)\|_{\infty}^{2}
                  \left(\sum_{|\alpha|=2} \|D^\alpha u_s(s)\|_2 \right)\: ds\\
   &      &  +CN_2(t)\left(\int_{t_0}^t \|u_s^{(p-1)/2}(s)Du_s(s)\|_2^2\: ds\right)^{1/2}.
\end{eqnarray*}
We use again (\ref{conv}) to derive
\begin{eqnarray*}
N_2(t) & \leq & 2C\sum_{|\alpha|=2}\|DD^\alpha u(t_0)\|_2\\
       &      & +2C\int_{t_0}^t\|u_s(s)\|_{\infty}^{2}
                  \left(\sum_{|\alpha|=2} \|D^\alpha u_s(s)\|_2 \right)\: ds
\end{eqnarray*}
for sufficiently large $t\geq t_0.$ Let us introduce
\[
N_3(t)=\sum_{|\alpha|= 2}\sup_{s\in[t_0,t]}\|DD^\alpha
u(s)\|_2,\qquad t\geq t_0.
\]
Then we can write
\begin{eqnarray}
N_2(t) & \leq & 2C\sum_{|\alpha|=2}\|DD^\alpha u(t_0)\|_2+2CN_3(t)\left(\int_{t_0}^t\|Du(s)\|_{\infty}^2 \:
       ds\right). \label{N2(t)}
\end{eqnarray}

We will combine the above estimate with Lemma \ref{lem02} (a) for the interval $[t_0,t]$:
\begin{eqnarray*}
\|DD^{\alpha}u(t)\|_2 & \leq &  C\sum_{|\alpha|=2} \|DD^\alpha u(t_0)\|_2\nonumber\\
                      &      & +  C\int_{t_0}^t \|u_s^{p-1}(s)D^{\alpha}u_s(s)\|_2 \: ds
                                \nonumber\\
                      &      &+ C\int_{t_0}^t \|u_s^{p-2}(s)(Du_s(s))^2\|_2 \: ds,
                               \nonumber
\end{eqnarray*}
where $t\geq t_0.$ It follows from (\ref{uni_est1}) and (\ref{uni_est2}) that
\begin{eqnarray*}
\|DD^{\alpha}u(t)\|_2 & \leq &  C\sum_{|\alpha|=2} \|DD^\alpha u(t_0)\|_2\nonumber\\
                      &      & +  C\int_{t_0}^t \|Du(s)\|_{\infty}^2\|D^{\alpha}u_s(s)\|_2  ds
                                \nonumber\\
                      &      &+ C\left(\int_{t_0}^t \|Du_s(s)\|_{\infty}^2 ds\right)^{1/2}
                                     \left(\int_{t_0}^t \|u_s^{(p-1)/2}(s)Du_s(s)\|_2^2 ds\right)^{1/2}.
                               \nonumber
\end{eqnarray*}
Summing over $|\alpha|=2$ and applying Lemma~\ref{lem01},  we simplify the estimate to
\begin{eqnarray}
N_3(t) & \leq &  C\sum_{|\alpha|=2} \|DD^\alpha u(t_0)\|_2\nonumber\\
                      &      & +  CN_3(t)\left(\int_{t_0}^t \|Du(s)\|^2_\infty \: ds\right)
                                \nonumber\\
                      &      &+CN_2(t)\left(\sum_{|\alpha|=1}\|DD^\alpha u(0)\|_2\right).
                               \nonumber
\end{eqnarray}
Since the integral
\[
\int_{0}^\infty\|Du(s)\|_{\infty}^2 \: ds
\]
is convergent, we choose a sufficiently large $t_0$ to finally get
\begin{eqnarray}
N_3(t) & \leq &  C\sum_{|\alpha|=2} \|DD^\alpha u(t_0)\|_2
                                 +CN_2(t)\left(\sum_{|\alpha|=1}\|DD^\alpha u(0)\|_2\right)
                                 \label{N3(t)}
\end{eqnarray}
for $t\geq t_0.$ Estimates (\ref{N2(t)}), (\ref{N3(t)}) are sufficient to bound
uniformly $N_2(t)+N_3(t).$ \qed

\section{Global existence of radial solutions in $H^{k}\times H^{k-1}$ with $k>3$}

Here we consider equation (\ref{main}) with a supercritical integer $p=2m+1$, $m\in \N$.
For $k\geq 3$ and $t\in[0,T_*)$, let us define
\[
E_{k}(t)=\sum_{|\alpha|=k-1}\|DD^{\alpha}u(t)\|_2,\qquad
N_{k}(t)=\sum_{1\le|\beta|\le k-1}\left(\int_{0}^{t}\|D^{\beta}u(s)\|_{\infty}^2 ds\right)^{1/2},
\]
where $u$ is the local solution of (\ref{main}), (\ref{data}) constructed in Lemma \ref{lem00}.

{\bf Proof of Theorem \ref{th01}.} We use induction in $k$.
From Theorem \ref{th1}, we know that the following hold when $k=3$:
\begin{equation}
\label{ind}
T_*=\infty,\quad \sup_{t\in[0,\infty)}E_k(t)<\infty,\quad \sup_{t\in[0,\infty)}N_k(t)<\infty.
\end{equation}
Assuming these claims for $k$, we will verify them for $k+1$. The proof starts from applying $D^{\alpha}$ to (\ref{main})
with $p=2m+1$ and solving for the highest derivatives:
\begin{eqnarray}
\Box D^{\alpha} u &= &c_{m,1}u_t^{2m}D^{\alpha}u_t\nonumber\\
                  &  &+c_{m,2}u_t^{2m-1}\sum_{\alpha_1+\alpha_2=\alpha}D^{\alpha_1}u_tD^{\alpha_2}u_t \nonumber\\
                  &  &+c_{m,3}u_t^{2m-2}\sum_{\alpha_1+\alpha_2+\alpha_3=\alpha}D^{\alpha_1}u_tD^{\alpha_2}u_tD^{\alpha_3}u_t\label{main:k-th}\\
                   & & +\cdots\nonumber\\
& &+c_{m,l}u_t^{2m+1-l}\sum_{\alpha_1+\alpha_2+\cdots+\alpha_{l}=\alpha}D^{\alpha_1}u_t\cdots D^{\alpha_{l}}u_t=0, \nonumber
\end{eqnarray}
where $l=\min\{2m+1,k\}$, $c_{m,i}$ are constants and $|\alpha_i|\geq 1$ for $i=1,\cdots, l.$
Making use of Lemma \ref{lem02} (a), we obtain
\begin{eqnarray}
\|DD^{\alpha} u(t)\|_2 &\le& C\|DD^{\alpha}(0)\|_2 +C\int_{0}^{t} \|u_s^{2m}D^{\alpha}u_s\|_2ds\nonumber\\
&&+\sum_{\alpha_1+\alpha_2=\alpha}C\int_0^t \|u_s^{2m-1}D^{\alpha_1}u_sD^{\alpha_2}u_s\|_2ds\nonumber\\
&&+\sum_{\alpha_1+\alpha_2+\alpha_3=\alpha}C\int_0^t \|u_s^{2m-2}D^{\alpha_1}u_sD^{\alpha_2}u_sD^{\alpha_3}u_s\|_2ds\label{k-deriv1}\\
&& +\cdots\nonumber\\
&&+\sum_{\alpha_1+\alpha_2+\cdots+\alpha_{l}=\alpha}C\int_0^t \|u_s^{2m+1-l}D^{\alpha_1}u_s\cdots D^{\alpha_{l}}u_s\|_2\}ds.
\nonumber
\end{eqnarray}

Only the first integrand involves a derivatives of order $k+1$. From (\ref{H_Est2}), we get
\begin{equation}
\label{k-deriv2}
\|u_s^{2m}D^{\alpha}u_s\|_2\leq C\|u_s\|_{\infty}^2E_{k+1}(s).
\end{equation}

In the second integrand, there are two cases: $\max\{|\alpha_1|,|\alpha_2|\}=k-1$
and $\max\{|\alpha_1|,|\alpha_2|\}\leq k-2$. We rely on (\ref{H_Est2}) to derive
\begin{eqnarray}
\|u_s^{2m-1}D^{\alpha_1}u_sD^{\alpha_2}u_s\|_2
&\le& \|u_s\|_{\infty}^{2m-2}\|u_s\|_{\infty} \|Du_s\|_{\infty}E_k(s)\nonumber\\
&\le& C\left(\sum_{1\le|\beta|\le2}\|D^{\beta}u(s)\|_{\infty}^2\right)E_k(s)\label{k-deriv3}
\end{eqnarray}
and
\begin{eqnarray}
\|u_s^{2m-1}D^{\alpha_1}u_sD^{\alpha_2}u_s\|_2
&\le& \|u_s\|_{\infty}^{2m-2}\|u_s\|_{\infty}\|D^{\gamma}u_s\|_{\infty}E_{k-|\gamma|+1}(s)\nonumber\\
&\le& C\left(\sum_{1\le|\beta|\le k-1}\|D^{\beta}u(s)\|_{\infty}^2\right)E_{k-|\gamma|+1}(s),\label{k-deriv4}
\end{eqnarray}
respectively, where $\gamma$ is such that $|\gamma|=\max\{|\alpha_1|,|\alpha_2|\}$.
\par\noindent

The third integrand also admits $\max\{|\alpha_1|,|\alpha_2|,|\alpha_3|\}=k-2$.
Let $\gamma$ be the multiindex satisfying $|\gamma|=k-2.$ Applying (\ref{H_Est2}), we have
\begin{eqnarray}
\|u_s^{2m-2}D^{\alpha_1}u_sD^{\alpha_2}u_sD^{\alpha_3}u_s\|_2
&\le& \|u_s\|_{\infty}^{2m-2}\|D^{\gamma}u_s\|_{\infty}\|Du_s\|_{\infty}E_2(s)\nonumber\\
&\le&C\left(\sum_{1\le|\beta|\le k-1}\|D^{\beta}u(s)\|_{\infty}^2\right).\label{k-deriv5}
\end{eqnarray}

Finally, we consider all integrands with $\max\{|\alpha_{j}|: j=1,2,\ldots, l\}\leq k-3.$
Since $k\ge4$, we can use the Sobolev embedding and (\ref{H_Est2}) to obtain
\begin{eqnarray}
&&\|u_s^{2m+1-l}D^{\alpha_1}u_sD^{\alpha_2}u_s\cdots D^{\alpha_{j}}u_s\|_2\nonumber\\
&\le& \|u_s\|_{\infty}^{2m+1-l}\left(\sum_{1\le|\beta|\le k-2}\|D^{\beta}u(s)\|_{\infty}^2\right)
\left(\sum_{j=1}^{k}E_{j}(s)\right)^{l-2}\label{k-deriv6}
\\
&\le& C\left(\sum_{1\le|\beta|\le k-2}\|D^{\beta}u(s)\|_{\infty}^2\right).\nonumber
\end{eqnarray}

We combine the basic estimate (\ref{k-deriv1}) with estimates (\ref{k-deriv2})--(\ref{k-deriv6}). This yields
\begin{eqnarray*}
 \|DD^\alpha u(t)\|_2 & \leq & C\|DD^\alpha u(0)\|_2
+ C\int_0^t \|u_s(s)\|_{\infty}^2 E_{k+1}(s) + CN_{k}(t)+C_k
\end{eqnarray*}
for all $t\in [0,T_*)$, where the constant $C_k=C_k(\|u_0\|_{H^{k}},\|u_1\|_{H^{k-1}})$.
Adding all such estimates with $|\alpha|=k$ and using (\ref{ind}), we have
\begin{eqnarray*}
 E_{k+1}(t)\leq C_{k+1}+C\int_0^t \|u_s(s)\|_\infty^2 E_{k+1}(s) \: ds,
\end{eqnarray*}
for some constant $C_{k+1}=C_{k+1}(\|u_0\|_{H^{k+1}},\|u_1\|_{H^{k}}).$ From the Gronwall inequality,
\[
E_{k+1}(t)\leq C_{k+1}\exp\left(\int_{0}^{t}\|u_s(s)\|_\infty^2ds\right).
\]
This shows that $E_{k+1}(t)$ is uniformly bounded on every finite subinterval of
$[0,T_\ast).$ Hence, $u$ can be continued to a global solution, such that
$\sup_{t\in [0,\infty)}E_{k+1}(t)<\infty.$

It remains to check the uniform estimates in $L_t^2L_x^{\infty}$.
We apply Lemma \ref{lem02} (b) to the solution of (\ref{main:k-th}). The calculations are very similar to the above ones, so we
give only the final estimate:
\begin{eqnarray*}
  \sum_{|\alpha|=k}\left(\int_{0}^t \|D^\alpha u(s)\|_\infty^2\: ds\right)^{1/2}
   & \leq & C\int_0^t \|u_s(s)\|_\infty^2 E_{k+1}(s) \: ds+ C_{k+1}.
\end{eqnarray*}
Noticing that  $\sup_{t\in[0,\infty)}E_{k+1}(t)<\infty$ and $\sup_{t\in[0,\infty)}N_1(t)<\infty$, we obtain
the final estimate in (\ref{ind}) with $k+1$, i.e., $\sup_{t\in[0,\infty)}N_{k+1}(t)<\infty$.
The proof of higher regularity by induction is complete.

Similarly, we can show that $C^\infty$ regularity is preserved during the evolution of compactly supported radial data.
 \qed

\begin{center}
ACKNOWLEDGMENTS
\end{center}
The authors are very grateful to Professor Hideo Kubo for his useful comments.

\vskip 0.5truecm

% ----------------------------------------------------------------
%\bibliographystyle{siamplain}
%\bibliography{references}

\end{document}